\documentclass[12pt]{article}
\usepackage{amssymb,amsmath,amsthm,latexsym}

\usepackage{amsfonts}

\usepackage{amscd}
\usepackage{amsbsy}

\usepackage{color}
\usepackage[normalem]{ulem}

\title{Hyperovals of $H(3,q^2)$ when $q$ is even}
\author{Antonio Cossidente\\
Dipartimento di Matematica e Informatica\\ Universit\`a della Basilicata\\
I-85100 Potenza -- Italy\\\texttt{antonio.cossidente@unibas.it}
\\
\\
Oliver H. King\\School of Mathematics and Statistics\\
Newcastle University\\
Newcastle upon Tyne NE1 7RU\\
England, UK
\\{\tt
oli.king@ncl.ac.uk}
\\
\\
Giuseppe Marino\thanks{This work was
supported by the Research Project of MIUR (Italian Office for
University and Research) ``Geometrie su Campi di Galois, piani di traslazione e geometrie di incidenza'', by the Research group GNSAGA of INDAM and by a grant of the London Mathematical Society}\\Dipartimento di matematica\\Seconda Universit\`a di Napoli\\ I-81100 Caserta -- Italy\\\texttt{giuseppe.marino@unina2.it}}

\date{}

\newtheorem{theorem}{Theorem}[section]

\newtheorem{lemma}[theorem]{Lemma}

\newtheorem{remark}[theorem]{Remark}

\newtheorem{prop}[theorem]{Proposition}

\def\cC{\mathcal C}

\def\cH{\mathcal H}
\def\cI{\mathcal I}

\def\cL{\mathcal L}
\def\cM{\mathcal M}

\def\cO{\mathcal O}

\def\cQ{\mathcal Q}

\def\cW{\mathcal W}

\def\0{\mathbf 0}

\def\GammaU{{\rm \Gamma U}}
\def\Prf{\textbf{Proof. }}

\def\<{\langle}
\def\>{\rangle}
\newcommand{\Tr}{\mathrm{Tr}}

\begin{document}
\maketitle
\newpage
\abstract{For even $q$, a group $G$ isomorphic to $PSL(2,q)$
stabilizes a Baer conic inside a symplectic subquadrangle ${\cal
W}(3,q)$ of ${\cal H}(3,q^2)$. In this paper the action of $G$ on
points and lines of ${\cal H}(3,q^2)$ is investigated. A
construction is given of an infinite family of hyperovals of size
$2(q^3-q)$ of ${\cal H}(3,q^2)$, with each hyperoval having the
property that its automorphism group contains $G$. Finally it is
shown that the hyperovals constructed are not isomorphic to known
hyperovals.}

\bigskip


\section{Introduction}

A hyperoval (or {\it local subspace}) $\cO$ in a polar space $\cal P$ of rank $r\ge 2$ is defined to be a set of points of $\cal P$ with the property that each line of $\cal P$ meets $\cO$ in either $0$ or $2$ points.
Hyperovals of polar spaces are related to locally polar
spaces. From a result of Buekenhout and Hubaut
\cite[Proposition 3]{BH} it follows that if $A$ is a polar space of polar rank $r$,
$r\ge 3$ and order $n$, and $H$ a hyperoval of $A$ then $H$
equipped with the graph induced by $A$ on $H$, is the adjacency
graph of a locally polar space of polar rank $r-1$ and order $n$
such that the residual space $H_P$ at any point $P\in H$ is
isomorphic to Cone$_P(A)$. This result makes interesting the
classification of all local subspaces of polar spaces.

Recently, B. De Bruyn \cite{BdB} proved the following theorem.

    \begin{theorem}
Let ${\cal S}$ be a generalized quadrangle of order $(s,t)$ and
let $H$ be a hyperoval of ${\cal S}$. Then
\begin{itemize}
\item[1.] $2$ is a divisor of $\vert H\vert$; \item[2.] We have
$\vert H\vert\ge 2(t+1)$, with equality if and only if there
exists a regular pair $\{x,y\}$ of non--collinear points of $\cal
S$ such that $H=\{x,y\}^\perp\cup\{x,y\}^{\perp\perp}$; \item[3.]
We have $\vert H\vert\ge(t-s+2)(s+1)$. If equality holds then
every point outside $H$ is incident with exactly $(t-s)/2+1$ lines
which meet $H$; \item[4.] We have $\vert H\vert\le 2(st+1)$, with
equality if and only if $H$ is a $2$--ovoid.
\end{itemize}
    \end{theorem}

As observed in \cite[Remark 2, p. 404]{BH} when $r=2$ we can still say that
a hyperoval of a generalized quadrangle $\cal S$ is a graph of degree equal to $\vert$Cone$_P({\cal S})\vert$ which has
the property to be triangle free. More precisely, a subgraph of a point $GQ(s,t)$--graph is
a hyperoval if and only if it is a regular graph of valency $t+1$ with an even number of vertices and has no triangles.

In this paper we shall focus on the generalized quadrangle ${\cal
H}(3,q^2)$, $q$ even, the incidence structure of all points and
lines (generators) of a non-singular Hermitian surface in
$PG(3,q^2)$, a generalized quadrangle of order $(q^2,q)$, with
automorphism group $P\Gamma U(4,q^2)$, see \cite{BS} for more
details.

The hyperovals of ${\cal H}(3,4)$ have been classified by Makhnev
\cite{M}. He found hyperovals of size $6,8,10,12,14,16,18$.

For some information on hyperovals of Hermitian variety ${\cal
H}(n,q^2)$, $n=4,5$, see \cite{BdB}, \cite{CM} and \cite{CMP}.

In this paper we study the action of the stabilizer $G$ of a Baer conic inside a symplectic subquadrangle ${\cal W}$ of ${\cal H}(3,q^2)$, $q$ even, on points and lines of ${\cal H}(3,q^2)$, which is interesting in its own right, and construct an infinite family of
hyperovals of ${\cal H}(3,q^2)$, $q$ even, of size $2(q^3-q)$ admitting the group $G$ as an automorphism group. Finally, we investigate
the isomorphism problem with the known hyperovals on ${\cal H}(3,q^2)$.

    \section{The geometric setting}\label{prel}

For the reader's convenience, some facts about the Hermitian
surface are summarized below.

In $PG(3,q^2)$ a {\em non--singular Hermitian surface} is
defined to be the set of all absolute points of a non--degenerate
unitary polarity, and is denoted by ${\cal H}(3,q^2)$.

A Hermitian surface ${\cal H}(3,q^2)$ has the
following properties.
\begin{enumerate}
\item[1.] The number of points on ${\cal H}(3,q^2)$ is
$(q^2+1)(q^3+1)$. \item[2.] The number of generators of ${\cal
H}(3,q^2)$ is $(q+1)(q^3+1)$. \item [3.] Any line of $PG(3,q^2)$
meets $\cal H$ in either $1$ point ({\em tangent}) or $q+1$ points
({\em hyperbolic}) or $q^2+1$ points ({\em generator}). \item[4.]
Through every point $P$ of ${\cal H}(3,q^2)$ there pass exactly
$q+1$ generators, and these generators are coplanar. The plane
containing these generators, say $\pi_P$, is the polar plane of
$P$ with respect to the unitary polarity defining ${\cal
H}(3,q^2)$. The tangent lines through $P$ are precisely the
remaining $q^2-q$ lines of $\pi_P$ incident with $P$, and $\pi_P$
is called the {\em tangent plane} to ${\cal H}(3,q^2)$ at $P$.
\item[5.]Every plane of $PG(3,q^2)$ which is not a tangent plane
to ${\cal H}(3,q^2)$ meets ${\cal H}(3,q^2)$ in a non-degenerate
Hermitian curve ({\em secant plane}).
\end{enumerate}

We assume throughout that $q$ is {\bf even}. We start by setting up an appropriate co-ordinate system. We begin with a co-ordinate system for a 3--dimensional space $PG(3,q)$ containing an elliptic quadric $\cQ=\cQ^-(3,q)$ which in turn contains a non-degenerate conic $\cC$. The co-ordinate system may be chosen so that $\cQ$ has equation \[X_0X_3 + X_1^2 + \sigma X_1 X_2 +X_2^2=0,\]
where $\sigma$ is an element of $GF(q)$ such that the polynomial $x^2 + \sigma x +1$ is irreducible over $GF(q)$, and so that $\cC$ is the intersection of $\cQ$ with the plane $X_2=0$; the points of $\cC$ are  $(0,0,0,1)$  and $(1,s,0,s^2)$ with $s \in GF(q)$ and $\cC$ has nucleus $N=(0,1,0,0)$. Let $\cW=\cW(3,q)$ be the symplectic space consisting of all the points of $PG(3,q)$ together with the totally isotropic lines with respect to the alternating form on $PG(3,q)$:
\[H((X_0,X_1,X_2,X_3),(Y_0,Y_1,Y_2,Y_3))=X_0Y_3+X_3Y_0 +\sigma (X_1Y_2 + X_2 Y_1).\]
The Hermitian surface $\cH=\cH(3,q^2)$ is obtained by extending $H$ to a Hermitian form (also denoted by $H$) on $PG(3,q^2)$ given by:
\[H((X_0,X_1,X_2,X_3),(Y_0,Y_1,Y_2,Y_3))=X_0Y_3^q+X_3Y_0^q +\sigma (X_1Y_2^q + X_2 Y_1^q),\]
and the Hermitian surface then has equation
\[X_0X_3^q+X_3X_0^q +\sigma (X_1X_2^q + X_2 X_1^q)=0\]
and it contains all the points of $\cW$. We denote by $\tau$ the semilinear map: $(X_0,X_1,X_2,X_3) \to (X_0^q,X_1^q,X_2^q,X_3^q)$ arising from the Frobenius automorphism $X \to X^q$ of $GF(q^2)$. Then $\cW$ consists of all of the points of $PG(3,q^2)$ fixed by $\tau$.

Every generator of $\cH$ either meets $\cW$ in a totally isotropic
Baer subline of $\cW$ or it is disjoint from it. This means that a
point of $\cH\setminus \cW$ lies on a unique generator arising
from a totally isotropic Baer subline of $\cW$. Moreover, since
$q$ is even, the totally isotropic Bear sublines of $\cW$ are
tangents to the elliptic quadric $\cal Q$ (see \cite{BS}).

The group $P\GammaU(4,q^2)$ consists of all linear or semilinear
collineations fixing $\cH$ and contains the group $PSp(4,q)$
(isomorphic to $Sp(4,q)$) acting on $\cW$; in particular it
contains the mapping $\tau$. In turn the group $PSp(4,q)$ contains
$P\Omega^-(4,q)$ acting on $\cQ$. The stabilizer in
$P\Omega^-(4,q)$ of the conic $\cC$ is a group, denoted by $G$,
that is isomorphic to $PSL(2,q)$ (and so also to $SL(2,q)$). Note
that $G$ acts transitively on the points of $\cC$ and fixes the
nucleus $N$. We shall find it helpful to work with the elements of
$G$ as matrices in  $\Omega^-(4,q)$. We shall consider the points
as column vectors, with matrices in $G$ acting on the left.

\begin{prop}\label{2.1}
The elements of $G$ are given by
\[\left( \begin{matrix}
a^2 & 0 & \sigma ab & b^2 \\
ac & 1 & \sigma bc & bd \\
0 & 0 & 1 & 0 \\
c^2 & 0 & \sigma cd & d^2
 \end{matrix} \right)\]
where $a,b,c,d \in GF(q)$ and $ad+bc=1$.
\end{prop}
\Prf It is readily seen that $\left( \begin{matrix} a & b \\ c & d   \end{matrix} \right)  \mapsto  \left( \begin{matrix}
a^2 & 0 & \sigma ab & b^2 \\
ac & 1 & \sigma bc & bd \\
0 & 0 & 1 & 0 \\
c^2 & 0 & \sigma cd & d^2
 \end{matrix} \right)$ is an isomorphism from
 $SL(2,q)$ to a subgroup $G'$ of $GL(4,q)$, and a straightforward
calculation shows that the elements of $G'$ preserve the quadratic
form $Q(X_0,X_1,X_2,X_3)=X_0X_3 + X_1^2 + \sigma X_1 X_2 +X_2^2$
of $Q^-(3,q)$. In other words $G' \leq O^-(4,q)$. It is clear that
$G'$ fixes the plane $X_2=0$ and therefore it stabilizes $\cC$. If
$q>2$, then the intersection of $G'$ with $\Omega^-(4,q)$ is a
normal subgroup of $G'$ of index $1$ or $2$, but $G'$ is simple,
which implies that $G' \leq \Omega^-(4,q)$. If $q=2$, the $G'$ is
no longer simple but can still be shown to lie in $\Omega^-(4,q)$.
It follows that $G'=G$.

\subsection{Point and line orbits of $PSL(2,q)$ on $\cH(3,q^2)$}
In considering the orbits of $G$ on the set of lines of $\cH$, we note that each orbit consists either of lines that meet $PG(3,q)$ in  lines of $\cW$ or of lines that are skew to $\cW$.
\begin{prop}\label{2.2}
The group $G$ has the following orbits on lines of $\cH$:

  \begin{enumerate}
     \item An orbit of length $q+1$ consisting of tangent lines to the conic $\cC$ that pass through the nucleus $N$.
     \item $q$ orbits of length $q+1$, each consisting of tangent lines to the conic that do not pass through the nucleus $N$.
     \item An orbit of length $q^3-q$ consisting of lines that are tangent to $\cQ$ but do not meet $\cC$.
     \item $2q$ orbits of length $(q^3-q)/2$ lines, each consisting of lines that are disjoint from $\cW$.
  \end{enumerate}
\end{prop}
\Prf
 \begin{enumerate}
     \item We know that $G$ fixes $N$ and acts transtively on the $q+1$ points of $\cC$ and so it acts transitively on the $q+1$ tangents to $\cC$ that pass through $N$.
     \item Suppose that $m$ is a line of $\cH$ that is tangent to $\cC$ but does not pass through $N$. Then the orbit containing $m$ includes tangents to each point of $\cC$. In particular it includes a tangent to the point $(1,0,0,0)$: we assume that $m$ passes through this point. The tangent plane to $\cH$ at $(1,0,0,0)$ has equation $X_3=0$ and this plane contains the $q+1$ generators of $\cH$ that pass through $(1,0,0,0)$. The generator $m$ must pass through $(0,\lambda,1,0)$ for some $\lambda \in GF(q)$. Let $G_1$ be the stabilizer in $G$ of the point $(1,0,0,0)$, then $G_1$ consists of matrices of the form: $A=\left( \begin{matrix}
a^2 & 0 & \sigma ab & b^2 \\
0 & 1 & 0 &bd \\
0 & 0 & 1 & 0 \\
0 & 0 & 0 & d^2
 \end{matrix} \right)$ with $ad=1$. Under such an element of $G_1$, $(0,\lambda,1,0)$ is mapped to $(\sigma ab,\lambda,1,0)$, which still lies in $m$. It follows that $G_1$ is the stabilizer in $G$ of $m$ and the line orbit of $G$ containing $m$  contains exactly one line passing through each point of $\cC$.
     \item Note that there are $q^2-q$ points of $\cQ \setminus \cC$ and through each of these points there pass $q+1$ generators. Each generator of $\cW$ is tangent to $\cQ$. Hence there are $q^3-q$ lines that are tangent to $\cQ$ but do not meet $\cC$. Let us choose one of these lines, $m: X_1=0, X_0=X_3$ which passes through the point $(1,0,1,1)$, and consider its stabilizer in $G$. Suppose that $A=\left( \begin{matrix}
a^2 & 0 & \sigma ab & b^2 \\
ac & 1 & \sigma bc & bd \\
0 & 0 & 1 & 0 \\
c^2 & 0 & \sigma cd & d^2
 \end{matrix} \right)$ maps $(0,0,1,0)$ to a point of $m$ and fixes $(1,0,1,1)$. Then $bc=0$, implying that $ad=1$, and also $ab=cd=0$, from which it follows that $b=c=0$; it also follows that $a^2=d^2=1$, i.e., $a=d=1$ and $A=I_4$.
 The Orbit-Stabilizer Theorem tells us that  the $q^3-q$ lines that are tangent to $\cQ$ but do not meet $\cC$ form an orbit under
$G$.
     \item In the parts above we have identified all $q^3+q^2+q+1$ lines that meet $PG(3,q)$ and so there remain $q(q^3-q)$ lines, all of which are skew to $PG(3,q)$. The following argument is prompted by \cite{AC2}. Consider the  line-orbit $\cL_1$ of length $q+1$ consisting of tangent lines to the conic $\cC$ that pass through the nucleus $N$. Altogether these lines contain $(q+1)(q^2-q)$ points of $\cH \setminus \cW$. Through each of the points there pass $q+1$ lines of $\cH$, of which one is in $\cL_1$ and the remainding $q$ lines are skew to $\cW$. Let $\cM$ be the set of lines so constructed (excluding the lines of $\cW$). No line of $\cM$ can meet two lines of $\cL_1$, because $\cH$ does not contain triangles. Thus $\cM$ contains  $q(q+1)(q^2-q)$ lines, i.e., all the lines of $\cH$ skew to $\cW$. Each orbit of $G$ on $\cM$ contains a line $\ell$ that meets the line $\ell_1 \in \cL_1$ given by $X_2=X_3=0$. Let $R=(1, \xi,0,0)$  (with $\xi \notin GF(q)$) be the point of intersection. Then the stabilizer in $G$ of $\ell$ stabilizes $R$ and therefore stabilizes $\ell_1$ and also $(1,0,0,0)$. An element of $G$ that fixes $(1,0,0,0)$ is given by $A=\left( \begin{matrix}
a^2 & 0 & \sigma ab & b^2 \\
0 & 1 &0 & bd \\
0 & 0 & 1 & 0 \\
0 & 0 &0 & d^2
 \end{matrix} \right)$ with $ad=1$.  If such a matrix fixes $R$, then $a=d=1$, but $b$ can take any value in $GF(q)$ (so the stabilizer in $G$ of $R$ has order $q$). The line $\ell$ contains a point with co-ordinates $(\nu,0,\xi,\sigma \xi^{q+1})$ for some $\nu \in GF(q)$. The image of this point under multiplication by $A=\left( \begin{matrix}1 & 0 & \sigma b & b^2 \\
0 & 1 &0 & b \\
0 & 0 & 1 & 0 \\
0 & 0 &0 &1
 \end{matrix} \right)$ is $(\nu + \sigma b \xi + b^2 \sigma \xi^{q+1}, b \sigma \xi^{q+1}, \xi, \sigma \xi^{q+1})$, which lies on $\ell$ when $b \sigma \xi^{q+1}=\xi (\sigma b \xi + b^2  \sigma \xi^{q+1})$. The two solutions to this equation are $b=0$ and $b=(\xi + \xi^{q} )/\xi^{q+1} \neq 0$.  Hence the stabilizer in $G$ of $\ell$ has order $2$. It follows that $\ell$ lies in an orbit of length $(q^3-q)/2$.
  \end{enumerate}

\begin{prop}\label{2.3}
The group $G$ has the following orbits on points of $\cH \setminus \cW$:

  \begin{enumerate}
     \item $q$ orbits of length $q^2-1$ consisting of points lying on  tangent lines to the conic $\cC$ that pass through the nucleus $N$.
     \item  $q$ orbits of length $q^3-q$, each consisting of  points  lying on  tangent lines to the conic $\cC$ that do not pass through the nucleus $N$.
     \item $q^2-q$ orbits of length $q^3-q$, each containing exactly one point of each of the lines in the line-orbit of length $q^3-q$.
  \end{enumerate}
\end{prop}
\Prf Note that each line of $\cW$ is tangent to exactly one point of $\cQ$ and each point of $\cQ$ lies on $q+1$ tangents lines in $\cW$. Furthermore a point of $\cH \setminus \cW$ lies on exactly one `extended' line of $\cW$, i.e., a line of $\cH$ that contains a line of $\cW$.
 \begin{enumerate}
     \item Let $\cL_1$ and $\ell_1$ be as in part (d) of Proposition \ref{2.2}. Then we know that the stabilizer in $G$ of the point  $R=(1,\xi,0,0)$ (with $\xi \notin GF(q)$) has order $q$, which implies that this point lies in an orbit of length $q^2-1$. Given that there are $q^3-q$ points of $\cH \setminus \cW$ lying on lines in the $G$-orbit  of $\ell_1$, there are $q$ orbits of such points.

     \item  Suppose that $\cL_2$ is a line-orbit of length $q+1$ consisting of tangent lines to the conic $\cC$ that do not pass through the nucleus $N$ and consider a line $m \in \cL_2$ that passes through $(1,0,0,0)$. As we have seen in Proposition \ref{2.2} (b), such a line passes through $(0,\lambda,1,0)$ for some $\lambda \in GF(q)$, and the stabilizer in $G$ of $m$ is the subgroup denoted $G_1$. A point of $\cH \setminus \cW$ on $m$ has coordinates  $(\xi,\lambda,1,0)$ with $\xi \notin GF(q)$. If $A=\left( \begin{matrix}
a^2 & 0 & \sigma ab & b^2 \\
0 & 1 &0 & bd \\
0 & 0 & 1 & 0 \\
0 & 0 &0 & d^2
 \end{matrix} \right)$ is an element of $G_1$ fixing $(\xi,\lambda,1,0)$, then $a^2 \xi + \sigma ab=\xi$, which implies that $a=1$ and $b=0$. Thus the stabilizer in $G$ of $(\xi,\lambda,1,0)$ has order $1$, and this implies that this point lies in an orbit of length $q^3-q$.
     \item   Suppose that $\cL_3$ is the line-orbit of length $q^3-q$ consisting of lines that are tangent to $\cQ$ but do not meet $\cC$, and that $\ell_3 \in \cL_3$. If $X$ is any point of $\ell_3$ not in $\cW$, then the $q^3-q$ lines $g(\ell_3)$ ($g \in G$) are disjoint outside $\cW$ and so the point-orbit of $G$ containing $X$ contains a point of each line in $\cL_3$. The orbit length cannot be greater than $q^3-q$, so it is exactly $q^3-q$ and the orbit contains exactly one point of each line in $\cL_3$.
  \end{enumerate}

\section{The construction of the hyperovals}\label{construction}

Consider the point-orbits identified in Proposition \ref{2.3}.\\

 \begin{enumerate}
     \item $q$ orbits of length $q^2-1$ consisting of points of $\cH \setminus \cW$ lying on tangent lines to the conic $\cC$ that pass through the nucleus $N$.\\
Let $\cO_1$ be one of these orbits. Then each point of $\cO_1$ lies on one of $q+1$ tangent lines to $\cC$ passing through $N$, and $G$ acts transitively on these lines so that each meets $\cO_1$ in the same number of points. It follows that each of these tangent lines meets $\cO_1$ in $q-1$ points. Hence this orbit cannot form a hyperoval.
     \item  $q$ orbits of length $q^3-q$, each consisting of  points of $\cH \setminus \cW$ lying on  tangent lines to the conic $\cC$ that do not pass through the nucleus $N$.\\
Let $\cO_2$ be one of these orbits. Then each point of $\cO_2$ lies on one of $q+1$ tangent lines to $\cC$ not passing through $N$, and $G$ acts transitively on these lines so that each meets $\cO_2$ in the same number of points. It follows that each of these tangent lines meets $\cO_2$ in $q^2-q$ points. Hence this orbit cannot form a hyperoval when $q>2$.

     \item $q^2-q$ orbits of length $q^3-q$, each containing exactly one point of each of the lines in the line-orbit of length $q^3-q$.\\
We consider a line $\ell_3$ that passes through $(1,0,1,1) \in \cQ \setminus \cC$ and $(0,0,1,0)$; this line lies in the line-orbit of $G$ of length $q^3-q$. A point on this line that does not lie in $\cW$ has co-ordinates $(\xi,0,1,\xi)$ for some $\xi \in GF(q^2) \setminus GF(q)$. Let $\cI=\cI_{\xi}$ be the point-orbit of $G$ containing the point $(\xi,0,1,\xi)$. We show that $\cI$ is a set of type $(0,1,2)$ with respect to lines of $\cH$ and determine when $\cI \cup \cI^{\tau}$ is a hyperoval.
  \end{enumerate}

\begin{prop}
The set $\cI_{\xi}$ lies on the quadric ${\cal Q}_{\xi}$ of
$PG(3,q^2)$ with equation $X_0X_3+X_1^2+\sigma
X_1X_2+\xi^2X_2^2=0$, and it is of type $(0,1,2)$ with respect to
generators of $\cH$.
\end{prop}
\Prf  Observe that ${\cal Q}_{\xi}$ contains the point
$(\xi,0,1,\xi)$. A straightforward calculation shows that $G$
fixes $\cQ_{\xi}$ and therefore $\cI_{\xi}$ lies on $\cQ_{\xi}$.
Now, we consider the four types of line indicated by the orbits of
$G$ on lines of $\cH$. The lines in the first three types of orbit
are all extended lines of $\cW$: we know that a point in
$\cI_{\xi}$ lies on exactly one extended line of $\cW$, and that
is a line of the third type which contains exactly on point of
$\cI_{\xi}$, so lines of the first three types meet $\cI_{\xi}$ in
at most one point. If $\ell$ is a line of $\cH$ that meets
$\cI_{\xi}$ and is not an extended line of $\cW$, then $\ell$ is
disjoint from $\cW$ and lies in a line-orbit $\cL_4$ of $G$ of
length $(q^3-q)/2$. Given that each line of $\cL_4$ meets
$\cI_{\xi}$ in the same number of points, and that each point of
$\cI_{\xi}$ lies on a line in $\cL_4$, we see that $\ell$ meets
$\cI_{\xi}$ in $2$ points. Therefore lines of the fourth type are
either disjoint from $\cI_{\xi}$ or meet it in $2$ points. Hence
$\cI_{\xi}$ is of type $(0,1,2)$ with respect to lines of $\cH$.

\begin{prop}\label{3.2}
Given $\xi \in GF(q^2) \setminus GF(q)$, the following are equivalent:
  \begin{enumerate}
     \item The quadric $\cQ_{\xi}$ is hyperbolic.
     \item The quadric $\cQ_{\xi^q}$ is hyperbolic.
     \item $\Tr_{GF(q^2)/GF(2)}\left(\dfrac{\xi^2}{\sigma^2}\right) =0$.
     \item The polynomial $x^2 + \dfrac{\sigma}{\xi} x + 1$ is reducible over $GF(q^2)$.
     \item $\Tr_{GF(q)/GF(2)}\left(\dfrac{\xi^2 + \xi^{2q}}{\sigma^2}\right) =0$.
     \item The polynomial $x^2 + \sigma x + \xi^2 + \xi^{2q}$ is reducible over $GF(q)$.
  \end{enumerate}
\end{prop}
\Prf
The quadric $\cQ_{\xi}$ is hyperbolic precisely when the polynomial $x^2 + \sigma x + \xi^2$ is reducible over $GF(q^2)$. This is equivalent to $x^2 + \sigma x + \xi^{2q}$ being reducible over $GF(q^2)$ and hence to $\cQ_{\xi^q}$ being hyperbolic.\\
It is known that the polynomial  $x^2 +ax+b$ is reducible over $GF(q)$ precisely when $\Tr_{GF(q)/GF(2)} \left(\dfrac{b}{a^2}\right) =0$, see e.g. \cite[p. 8]{JWPH}. Thus $\cQ_{\xi}$ is hyperbolic precisely when $\Tr_{GF(q^2)/GF(2)}\left(\dfrac{\xi^2}{\sigma^2}\right) =0$ and this is also the condition that $x^2 + \dfrac{\sigma}{\xi} x + 1$ is reducible over $GF(q^2)$.\\
The statement
$\Tr_{GF(q^2)/GF(2)}\left(\dfrac{\xi^2}{\sigma^2}\right) =0$ is
equivalent to \\$\Tr_{GF(q)/GF(2)}\left(\dfrac{\xi^2 +
\xi^{2q}}{\sigma^2}\right) =0$, which is precisely the condition
for $x^2 + \sigma x + \xi^2 + \xi^{2q}$ to be reducible over
$GF(q)$.

\begin{remark}{\rm
Notice that $\Tr_{GF(q^2)/GF(2)} x=0$ for all $x \in GF(q)$ and
for $q^2/2$ values of $x \in GF(q^2)$. Hence $\cQ_{\xi}$ is
hyperbolic for $(q^2-2q)/2$ values of $\xi \in GF(q^2) \setminus
GF(q)$ and elliptic for $q^2/2$ values of $\xi \in GF(q^2)
\setminus GF(q)$.}
\end{remark}

\begin{prop}
\begin{itemize}
\item If ${\cal Q}_\xi$ is an elliptic quadric, then $\cI_{\xi} \cup\cI_{\xi}^{\tau}$ is a hyperoval of $\cH$ of size $2(q^3-q)$.
\item If ${\cal Q}_\xi$ is a hyperbolic quadric, then $\cI_{\xi} \cup\cI_{\xi}^{\tau}$ is a set of type $(0,2,4)$ with respect to the lines of $\cH$.
\end{itemize}
\end{prop}
\Prf We observe that $\cI_{\xi}^{\tau}$ is the point-orbit of $G$ containing the point $(\xi^q,0,1,\xi^q)$ because the matrices in $G$ have entries in $GF(q)$. The point $(\xi^q,0,1,\xi^q)$ does not lie in $\cI_{\xi}$, so $\cI_{\xi} \cup\cI_{\xi}^{\tau}$ is a union of two equally-sized orbits of $G$ and therefore has
size $2(q^3-q)$.

We consider again the four types of line indicated by the orbits of $G$ on lines of $\cH$. As we have seen, the lines in the first two types do not meet $\cI_{\xi}$ and the same applies to $\cI_{\xi}^{\tau}$. Lines of the third type lie in a single orbit, and the representative given by $X_1=0, X_0=X_3$ (earlier denoted $\ell_3$) meets $\cI_{\xi}$ in just the point $(\xi,0,1,\xi)$ and $\cI_{\xi}^{\tau}$ in just the point $(\xi^q,0,1,\xi^q)$, so meets $\cI_{\xi} \cup\cI_{\xi}^{\tau}$ in $2$ points.

Now suppose that $\ell$ is a line of the fourth type meeting
$\cI_{\xi}$ in at least one point. We may suppose that $\ell$
passes through $(\xi,0,1,\xi)$ (because if $\cL_4$ is the
line-orbit of $G$ containing $\ell$, then each line in $\cL_4$
meets each of $\cI_{\xi}$ and $\cI_{\xi}^{\tau}$ in the same
number of points and so there is a line in $\cL_4$ passing through
$(\xi,0,1,\xi)$). Let $\cQ_{\xi}$ be the quadric in $PG(3,q^2)$
with equation $X_0X_3 + X_1^2 + \sigma X_1X_2 + \xi^2 X_2^2=0$
defined in the previous proposition.

Observe that the quadric $\cQ_{\xi^q}$ contains the set $\cI_{\xi}^{\tau}$. Let $\pi$ be the tangent plane to $\cH$ at $(\xi,0,1,\xi)$: it has equation $\xi^q X_0+ \xi^q X_3 + \sigma X_1=0$. Each generator of $\cH$ that passes through $(\xi,0,1,\xi)$ lies in the plane $\pi$. One of these is the line: $X_1=0, X_0=X_3$ which extends a line of $\cW$. Given that $(\xi,0,1,\xi)$ lies on exactly one such line, the other $q$ generators through $(\xi,0,1,\xi)$ (of which $\ell$ is one) are disjoint from $\cW$ and therefore meet $\cI_{\xi}$ in $0$ or $2$ points; in this case they each meet $\cI_{\xi}$ in $2$ points. Note that the polar form associated with $\cQ_{\xi^q}$ is given by
\[F((X_0,X_1,X_2,X_3),(Y_0,Y_1,Y_2,Y_3))=X_0Y_3+X_3Y_0 +\sigma (X_1Y_2 + X_2 Y_1),\]
so that the tangent plane to $\cQ_{\xi^q}$ at $(\xi^q,0,1,\xi^q)$ is precisely the plane $\pi$.

If ${\cal Q}_{\xi}$ and $\cQ_{\xi^q}$ are elliptic, then $\pi$ meets $\cQ_{\xi^q}$ in just the point $(\xi^q,0,1,\xi^q)$. It follows that $\ell$ is disjoint from $\cI_{\xi}^{\tau}$. Hence $\ell$ meets $\cI_{\xi} \cup\cI_{\xi}^{\tau}$ in $2$ points. Finally, a line of the fourth type that is disjoint from $\cI_{\xi}$ meets $\cI_{\xi}^{\tau}$ in $0$ or $2$ points, so meets $\cI_{\xi} \cup\cI_{\xi}^{\tau}$ in $0$ or $2$ points. Hence, in this case
the set $\cI_{\xi} \cup\cI_{\xi}^{\tau}$ is a hyperoval.

If ${\cal Q}_{\xi}$ and $\cQ_{\xi^q}$ are  hyperbolic, then  $\pi$ meets $\cQ_{\xi^q}$ in two lines through the point $(\xi^q,0,1,\xi^q)$. Consider the points $(\xi^q b^2,\sigma, 1, \xi^q/b^2)$ with $0 \neq b \in GF(q)$. They lie in ${\cal I}^\tau$ because they are the images of $(\xi^q,0,1,\xi^q)$ under elements of $G$ represented by matrices with $a=d=0$ and $bc=1$ (see Proposition \ref{2.1}). Such points lie in $\pi$ precisely when $\xi^{2q}b^4 +\sigma^2b^2 + \xi^{2q}=0$, i.e., when $b^2 + \dfrac{\sigma}{\xi^q} b + 1=0$. It follows from Proposition \ref{3.2} that $(\xi^q b^2,\sigma, 1, \xi^q/b^2)$ lies on $\pi \cap \cI_{\xi^q}$ for some $0 \neq b \in GF(q)$. Hence, the generator joining this point and $(\xi,0,1,\xi)$ meets ${\cal I}^\tau$ in four points.
Hence, in this case the set $\cI_{\xi} \cup\cI_{\xi}^{\tau}$ is of type $(0,2,4)$ with respect to generators of $\cH$.

    \section{The isomorphism problem}

To date, the known infinite families of hyperovals on the Hermitain surface are:

\begin{itemize}

\item

{\it Del Fra-Ghinelli-Payne's family}: a member in this family has size $2(q^3-q)$ and it is the union of two ovoids of ${\cal H}(3,q^2)$ sharing a chord;

\item

{\it Symplectic hyperovals:} Any $2$--ovoid of ${\cal W}(3,q)$
embedded in $\cH(3,q^2)$  gives rise to a hyperoval of size
$2(q^2+1)$ of ${\cal H}(3,q^2)$. Such hyperovals have been called
{\em symplectic hyperovals} in \cite{AC1}, and examples of
2--ovoids of ${\cal W}(3,q)$ have been constructed in \cite{BLP}
and \cite{CCEM}.

    \end{itemize}

A hyperoval on ${\cal H}(3,q^2)$ constructed in the previous
proposition has size $2(q^3-q)$ and hence it cannot be a
symplectic hyperoval. On the other hand, it has the same size of
the hyperovals described by Del Fra, Ghinelli and Payne. Of
course, the natural question is concerned with the existence of an
isomorphism between our hyperovals and the Del Fra-Ghinelli-Payne
hyperovals. Notice that the construction in \cite[Theorem
6.1]{DGP} involves the union of any two ovoids of ${\cal
H}(3,q^2)$ with a common chord deleted. Examples of such
hyperovals arise from the union of two non--degenerate plane
sections of ${\cal H}(3,q^2)$ (see \cite[Theorem 6.6]{DGP}) and
they obviously contain chords and hence $q+1$ collinear points. By
construction, our hyperovals lie on the union of two elliptic
quadrics of $PG(3,q^2)$ and hence, if $q>2$ and cannot be
isomorphic to any Del Fra-Ghinelli-Payne's hyperoval. When $q=2$,
MAGMA computations (\cite{magma}) show that our hyperovals do not
lie on the union of two secant planes of ${\cal H}(3,4)$.

    \begin{remark}
{\rm When $q=2$, each orbit ${\cal O}_2$ (see point (b) of Section
\ref{construction}) turns out to be a hyperoval of size $6$.}
    \end{remark}

\section{The automorphism group}
Our aim in this section is to determine the automorphism group of
a hyperoval $\cO=\cI_{\xi} \cup\cI_{\xi}^{\tau}$ in $P\Gamma
U(4,q^2)$, with $\xi$ necessarily chosen so that $\cQ_{\xi}$ is
elliptic. We shall see that it has the structure $G \times C_2
\times C_2$, being the product of $G$, the symmetry group on $N$
and the  field automorphism group $\{1,\tau\}$. Since this group
is the same for each relevant choice of $\xi$, the hyperovals
constructed are (in a sense) equivalent.

We denote by $\Delta$ the tangent plane to $\cH$ at $N$, and its
intersection with $PG(3,q)$ by $\Delta_0$. Thus $\Delta$ has
equation $X_2=0$ and $\Delta_0$ is a plane in $PG(3,q)$ containing
$\cC$ and $N$. Recall that each point of $\cH$ lies on $q+1$
generators of $\cH$: for a point of $\cW$ all of these generators
lie in $\cW$ and are tangent to $\cQ$, and for a point of $\cH
\setminus \cW$ one generator is an extended line of $\cW$ and the
remainder are skew to $\cW$. We make the following observations:

 \begin{enumerate}
    \item A point of $\cC$ lies on $q+1$ tangents to $\cC$.
    \item The point $N$ lies on $q+1$ tangents to $\cC$.
    \item A point of $\cQ \setminus \cC$ lies on $q+1$ tangents to $\cQ \setminus \cC$ and therefore on no tangents to $\cC$.
    \item If $R \in \Delta_0 \setminus (\cC \cup \{N\})$  and if $m$ is a generator through $R$, then $m$ lies in $\Delta$ if and only if $N$ lies on $m$. Thus there is one generator through $R$ that lies in $\Delta$, and this is the line $RN$ that is tangent to $\cC$. The remaining $q$ generators must be tangent to $\cQ \setminus \cC$.
    \item If $R \in PG(3,q) \setminus \Delta_0$ and if $R \notin \cQ$, then the tangent plane to $\cH$ at $R$ meets $\Delta_0$ at a line $\ell_0$ of $\Delta_0$ that is either external or secant to $\cC$. If $\ell_0$ is external to $\cC$, then $R$ lies on no tangents to $\cC$. If $\ell_0$ is secant to $\cC$, then $R$ lies on $2$ tangents to $\cC$, and the line $RN$ of $PG(3,q)$ contains $q$ points (all the points except $N$) that each lie on $2$ tangents to $\cC$.
    \item If $R \in \cH \setminus \cW$, then $R$ lies on a unique tangent to a point on $\cQ$ that might or might not lie on $\cC$, so $R$ lies on $0$ or $1$ tangents to $\cC$.
\end{enumerate}

\begin{lemma}
Let $H$ be the stabilizer  in $P\Gamma U(4,q^2)$ of $\cO=\cI_{\xi} \cup\cI_{\xi}^{\tau}$. Then $H \leq \< PSp(4,q), \tau \>$ and $H$ fixes $\cC$.
\end{lemma}
\Prf We observe that $\cO$ has $2(q^3-q)$ points, through each of which pass $q+1$ generators. Each of these generators meets $\cO$ in $2$ points. Hence there are $(q^3-q)(q+1)$ generators that meet $\cO$ in $2$ points. The remaining $(q+1)^2$ generators are skew to $\cO$. Since each point of $\cO$ lies on a unique tangent to $\cQ$, and  that tangent does not meet $\cC$, we conclude that the remaining  $(q+1)^2$ generators are the tangent lines to $\cC$. Thus $H$ fixes (globally) the set of tangents to $\cC$. It follows that for points of $\cH$, the number of tangents that pass through the point is an invariant under $H$.

The only points that lie on $q+1$ tangents to $\cC$ are the points of $\cC$ together with $N$. Moreover $N$ is collinear with all points of $\cC$, whereas each point of $\cC$ is collinear with $N$ and with no other point of $\cC$. Hence $H$ fixes $\cC$ and $N$, and also the plane $\Delta$.

The points of $\Delta_0 \setminus \{\cC \cup \{N\}\}$ each lie on one tangent to $\cC$ through $N$ and $q/2$ secants to $\cC$, whereas points of $\Delta \setminus \Delta_0$ lie on exactly one (extended) line of $\Delta_0$. Hence $H$ fixes $\Delta_0$.

Consider the set of points $R \in PG(3,q) \setminus \Delta_0$ such
that the line $\ell_0$ of $\Delta_0$ in (e) above is secant to
$\cC$. This is the set of points of $\cH$ that lie on $2$ tangents
to $\cC$, so it is fixed (globally) by $H$. In particular an
element of $H$ maps the line $RN$ of $PG(3,q)$ to a line of
$PG(3,q)$ through $N$. This is sufficient to establish that $H$
fixes $\cW$, i.e., $H \leq \< PSp(4,q), \tau \>$.

\begin{lemma}
The stabilizer in $\< PSp(4,q), \tau \>$ of $\cC$ is $G \times T \times \{1,\tau\}$, where $T$ is the group of symplectic transvections centred on $N$.
\end{lemma}
\Prf We begin by noting that $\tau$ fixes $\cC$ and commutes with $PSp(4,q)$, so it suffices to consider the stabilizer of $\cC$ in $PSp(4,q)$. This group acts on $\Delta_0$ as the orthogonal group $PO(3,q)$, in other words it acts as $G$. Therefore, for any $h \in Stab_{PSp(4,q)}\cC$, there exists $g \in G$ such that $gh$ acts as the identity on $\Delta_0$. Thus we now assume that $h$ acts as the identity on $\Delta_0$. Therefore  $h=\left( \begin{matrix}
\mu & 0 & u &  0\\
0 & \mu & v &  0\\
0 & 0 & w & 0 \\
0 & 0 & z & \mu
 \end{matrix} \right)$ for some $u,v,w,z, \mu \in GF(q)$ with $w, \mu \neq 0$. For such a matrix to lie in $PSp(4,q)$ it is necessary that $u=z=0$, $\mu=1$ and $w=1$. Thus $h$ is a symplectic transvection centred on $N$. Finally we note that all symplectic transvections centred on $N$ fix $\cC$ and commute with $G$.

\begin{theorem}
Suppose that $\xi \in GF(q^2) \setminus GF(q)$ such that $\cQ_{\xi}$ is elliptic. Then the stabilizer in $P\Gamma U(4,q^2)$ of $\cO=\cI_{\xi} \cup\cI_{\xi}^{\tau}$ is $G \times S \times \{1,\tau\}$, where $S$ is the symmetry group on $N$.
\end{theorem}
\Prf Given that $G$ and $\{1,\tau\}$ stabilize $\cO$, we need only determine which transvections centred on $N$ fix $\cO$. Let $t=\left( \begin{matrix}
1 & 0 & 0 &  0\\
0 & 1 & v &  0\\
0 & 0 & 1 & 0 \\
0 & 0 & 0 & 1
 \end{matrix} \right)$ and consider $t(P)$, where $P=(\xi,0,1,\xi) \in \cI_{\xi}$. If $t(P) \in \cI_{\xi} \cup \cI_{\xi}^{\tau}$, then $t(P) \in \cQ_{\xi} \cup \cQ_{\xi^q}$ and this last property is easy to check. Note that $t(P)=(\xi,v,1,\xi)$.

If $t(P) \in \cQ_{\xi}$, then $\xi^2 +v^2 +\sigma v + \xi^2=0$, i.e., $v=0$ or $\sigma$. If $v=\sigma$, then  $t$ is the symmetry on $N$, which we now denote by $s$. Observe that $s$ is the only non-trivial transvection that lies in $PO^-(4,q)$, and $\{1,s\}$ is a subgroup of $PO^-(4,q)$ that we refer to as the symmetry group centred on $N$.

By Proposition \ref{3.2}, the polynomial $x^2 +\sigma x + \xi^2 + \xi^{2q}$ is irreducible over $GF(q)$ and hence  $v^2 +\sigma v + \xi^2 + \xi^{2q} \neq 0$ for all $v \in GF(q)$, i.e., $t(P) \notin \cQ_{\xi^q}$.

It remains to check that $s$ stabilizes $\cO$. We see that $s(P)=g(P)$, where $g$ is the element of $G$ given by $\left( \begin{matrix}
0 & 0 & 0 &  1\\
0 & 1 & \sigma &  0\\
0 & 0 & 1 & 0 \\
1 & 0 & 0 & 0
 \end{matrix} \right)$. Given that $s$ commutes with both $G$ and $\tau$, it follows that $s$ stabilizes $\cO$.

\end{document}